\documentclass[12pt]{article}
\usepackage{amsmath,amssymb, plain, theorem, epsfig, graphicx,fancyhdr}
\usepackage{amsfonts, latexsym, multirow}
\newcommand{\var}{\mbox{\sc var}}

\newtheorem{thm}{\sc THEOREM}
\newtheorem{lem}{\sc LEMMA}

\begin{document}
\begin{center}
{\Large\bf A Partial Order on Uncertainty and Information }
\vspace{0.1in}

{\sc By} Jiahua Chen\\
Department of Statistics, 
University of British Columbia\\
 Vancouver, British Columbia, Canada V6T 1Z2\\
jhchen@stat.ubc.ca
\vspace{0.1in}
\end{center}
\centerline{{\sc Summary}}

Information and uncertainty are closely related and extensively
studied concepts in a number of scientific disciplines such as
communication theory, probability theory, and statistics. 
Increasing the information arguably reduces the uncertainty
on a given random subject.
Consider the uncertainty measure as the variance of a random
variable. Given the information that its outcome is in an interval, 
the uncertainty is expected to reduce when the interval shrinks.
This proposition is not generally true. In this paper, we provide a necessary
and sufficient condition for this proposition when the random
variable is absolutely continuous or integer valued. We also
give a similar result on Shannon information. 

\vspace{2ex}

\vspace{2ex}

\noindent
{\small \noindent {\rm Keywords and phrases:} 
Cauchy distribution, 
Conditional variance,
Dispersion order,
Gini's mean difference,
Likelihood ratio order,
Log-concavity, Normal distribution,
Shannon information, 
Stochastic order, Total positivity.

\vspace{2ex}

\section{Introduction}
Information and uncertainty and their relationship
are familiar notions in daily life, but quantifying them
is not easy. Probability theory provides a
platform for the study of random objects but does not provide
a universal index. Harris (1982) proposes the 
(relative) entropy of a probability distribution as such
an index. The notion of entropy goes
back at least as far as Shannon (1948) when uncertainty
and information were seen as identical: the
quantitative uncertainty $U(X)$ about a randomly
distributed object $X$ was thought of as the
amount of information that observing $X$ would provide,
since then all uncertainty about it would vanish.

If $X$ is a random variable or a random vector, then 
in statistical science its variance 
or covariance matrix $\var(X)$ is often regarded as an uncertainty 
measure.
The variance is simple and popular, and it is routinely
used to index the uncertainty in estimators.
For its viability as an uncertainty measure, it is natural
to ask whether increased knowledge about $X$ reduces
the uncertainty as measured in terms of the variance. 
In this vein, Zidek and van Eeden (2003) and
Chen, van Eeden, and Zidek (2010) show that 
the conditional variance $\var(X|~ |X|<x)$
is an increasing function of $x$ when $X$ has a normal distribution,
regardless of its mean and variance.
In fact, this result may be regarded as a straight consequence
of earlier results, as will be detailed.
Yet it is easy to find counterexamples where this conditional
variance of $X$ does not increase with $x$. The size of
the family of distributions for which this
monotonicity holds remained an unsolved problem.

In the search for an answer, the results in Burdett (1996) 
provide additional insight. For any given $X$, denote the conditional mean
$\mu(x) = E(X| X \leq x)$ and the conditional variance $\sigma^2(x) = \var(X | X \leq  x)$.
Both $\mu(x)$ and $\sigma^2(x)$ 
play important roles in economics, actuarial science, reliability
theory, and many other disciplines. 
Burdett (1996) provides a necessary and sufficient condition for
$\sigma^2(x)$ to be a monotonic function when $X$ is absolutely continuous
with finite mean and variance.
In particular, log-concavity of the density function
or of the cumulative distribution function of $X$ is sufficient, when
the second moment of $X$ is finite.
In a paper on dispersion orders, Mailhot (1987) shows that
when the cumulative distribution function of $X$ is log-concave, 
the conditional distribution of $X$ given $X < x$ has increasing dispersion 
order in $x$. This order implies monotonicity of $\sigma^2(x)$,
but this result is obtained under a stronger condition than that of Burdett (1996). 
Mailhot (1987) also contains a result that covers
the normal result of Chen, van Eeden, and Zidek (2010).
There are undoubtedly more papers that contain results implying the
monotonicity of $\sigma^2(x)$ under various conditions. 

In this paper, we advance the monotonicity of
the uncertainty and information measures on several fronts. 
In one respect, we improve on Mailhot (1987), Burdett (1996),
and Chen, van Eeden, and Zidek (2010)
by establishing a partial order on the conditional variance.
Let $\sigma^2(A) =   \var(X | X \in A)$ for any measurable set $A$.
We give a necessary and sufficient condition on $X$ with
an absolute continuous distribution
under which $\sigma^2(A) \leq \sigma^2(B)$ for any two intervals $A \subset B$.
We say that $\sigma^2(A)$ with this property is {\it partially monotonic}.
When $A$ is a finite interval, the conditional variance is always
well defined. Hence, this result is very general.
For some distributions such as Cauchy, we
establish partial monotonicity for special interval classes of $A$.
This result is particularly interesting because the
variance of the Cauchy distribution does not exist.
This result is presented in Section 2. The result is further
applied to integer-valued random variables.

A scientific proposition stands only if the result can be repeated independently.
In this view, the uncertainty in $X$ may be measured by the difference
in the outcomes from two independently conducted 
experiments under identical conditions. 
Let $X_1$ and $X_2$ be two such outcomes. The best uncertainty
measure might be a function of $X_1 - X_2$.
Additional information in the form of $X_1, X_2 \in A$ should 
reduce the uncertainty in general. Let $\varphi(u)$
be any increasing function in $|u|$.
We show that a sufficient condition for $E \{\varphi(X_1 - X_2)|X_1, X_2 \in A \}$
to be partially monotonic is that the density function of $X$ is
log-concave. Based on this result, we further show that the
conditional Shannon information of $X_1 - X_2$ is partially monotonic
under the same condition. Hence,
Shannon information based on $X_1 - X_2$ is 
another sensible information measure. We present these results in
Section 3. The paper ends with a short discussion in Section 4.


\section{Partial monotonicity of the conditional variance}
\label{sec2}

Let $X$ be an absolutely continuous random variable.
Denote its cumulative distribution
function by $F(x)$ and its density function as $f(x)$. 
We first investigate the conditions under which
$\var(X | 0 < X < b)$ is an increasing function of $b$. Based on this
result, we give a necessary and sufficient condition
for the partial monotonicity of  $\var(X| X \in A)$.

For a fixed value of $b > 0$, and for $x \in [0, b]$,
define $F_1(x) = \int_0^x \{F(t)- F(0)\} dt$  and
\[
F_2 (x) = \int_0^x F_1(t) dt 
= \int_{0 \leq s \leq t \leq x} \{F(s)-F(0)\} ds dt.
\]
Note that $F_1(x)$ is a specific antiderivative of $F(x)$,
and $F_2(x)$ is a specific antiderivative of $F_1(x)$.

Without loss of generality, we assume $F(0) = 0$
and $F(b) > 0$ in the following derivations.
Applying the technique of integration by parts, the
conditional mean
\[
\mu(b)
= E\{ X | 0 \leq X \leq b\} 
= \int_0^b x dF(x) /F(b)
= b - F_1(b)/F(b).
\]
Similarly, we find
\begin{eqnarray*}
 E\{ X^2 | 0 \leq X \leq b\} 
&=&
  \int_0^b x^2 dF(x) /F(b) \\
&=&
b^2 - 2 \int_0^b x dF_1(x)/F(b)\\
&=&
b^2 - 2 bF_1(b)/F(b) + 2 F_2(b)/F(b).
\end{eqnarray*}
Consequently, we have the following expression for the
conditional variance:
\[
\sigma^2(b) = \var \{X | 0 \leq X \leq b\}
=
\frac{2 F_2(b)}{F(b)} - \frac{F_1^2(b)}{F^2(b)}.
\]
Its derivative with respect to $b$ is given by
\[
\frac{f(b)}{F^3(b)}\{F_1^2(b) - F(b)F_2(b)\}.
\]
Hence, $\var \{X | 0 \leq X \leq b\}$ is an increasing function of $b$
if and only if $F_1^2(b) - F(b)F_2(b) \geq 0$ for all $b>0$.
This is equivalent to $F_2(b)$ being log-concave.
The above proof has closely followed that of Burdett (1996).

We now summarize the above derivation by a theorem in which
$F(0) = 0$ is no longer assumed.

\begin{thm}
\label{thm1}
Let $X$ be an absolutely continuous random variable with 
cumulative distribution function $F(x)$.
The conditional variance $\var \{X | 0 \leq X \leq b\}$ is an
increasing function of $b$ if and only if
\[
\int_{0 \leq x \leq y \leq b} \{F(x)-F(0)\} dx dy
\]
is log-concave.
\end{thm}

The above theorem easily generalizes from
conditioning on $0 < X < b$
to conditioning on $a < X < b$.
The generalization leads to partial
monotonicity of $\var(X| X \in A)$. The proof of the next theorem
is straightforward and omitted.

\begin{thm}
\label{thm2}
Let $X$ be an absolutely continuous random variable with 
cumulative distribution function $F(x)$.
The conditional variance
$\var \{X | a \leq X \leq b\}$ is
increasing in $b$ if and only if
\begin{equation}
\label{ConditionA}
\int_{a \leq x \leq y \leq b} \{F(x)-F(a)\} dx dy
\end{equation}
is log-concave in $b$, and it is decreasing in $a$
if and only if
\begin{equation}
\label{ConditionB}
\int_{a \leq x \leq y \leq b} \{F(b)-F(x)\} dx dy
\end{equation}
is log-concave in $a$.

When both conditions are satisfied for all $a, b \in C$ for
some convex set $C$, then
$\var \{X | X \in A\}$ is partially monotonic in interval $A$ such that
$A \subset C$.
\end{thm}

Many closely related results have been established
in the literature. Under nearly identical conditions, Burdett (1996) establishes the
mononicity of the conditional variance when $A$ takes the form $(- \infty, b]$. 
Theorem \ref{thm2} is more general because it is applicable to
situations where the mean and variance of $X$ do not exist, and 
to both finite and infinite intervals.  

Let $F(x)$ be a cumulative distribution function and let
$F^{-1}(\alpha) = \inf \{x: F(x) \geq \alpha\}$ for $0 < \alpha < 1$.
$F$ has higher dispersion order than $G$
if for any $0 < \alpha \leq \beta < 1$
\[
F^{-1}(\beta) - F^{-1}(\alpha) \geq G^{-1}(\beta) - G^{-1}(\alpha).
\]
When $F(x)$ is log-concave, Mailhot (1987) shows
that the conditional distribution of $X$ given $X < a$ 
has increased dispersion order in $a$ and hence increased
conditional variance.
Hence, Mailhot (1987)'s result implies the result of Burdett (1996),
though the latter provides a necessary and sufficient condition.
Mailhot (1987) further shows that if the density function $f$ is log-concave,
the conditional distribution of $X$ given $a < X < b$ decreases in $a$
and increases in $b$ in dispersion order. This result
implies the partial monotonicity presented in Theorem \ref{thm2},
and completely covers the normal result of Chen, van Eeden, and Zidek (2010).
Theorem \ref{thm2}, however, succeeds at giving a necessary and sufficient
condition. Theorem \ref{thm2} might be useful for establishing a
dispersion order more broadly in reverse.

It is natural to examine under what distributions Conditions (\ref{ConditionA})
and (\ref{ConditionB}) are satisfied. We now give a sufficient condition and
a few examples. The notion of log-concavity plays a key role.
Log-concavity has been studied thoroughly in mathematics.
The following is a well-known fact; for its proof, see Bagnoli and Bergstrom (2005).

\begin{lem}
\label{lem3}
If a function $f(x)$ is log-concave for $x \in (a, b)$, 
then the following antiderivative
\[
F(x) = \int_a^x f(t) dt
\]
is also log-concave for $x \in (a, b)$ whenever it is well defined.
\end{lem}

Note that $a = -\infty$ and/or $b = \infty$ are
special cases. According to this lemma, if a density function
is log-concave, so is its cumulative distribution function.
Clearly, if $F(x)$ is log-concave, so is $F(\alpha x +\beta)$ 
for any real numbers $\alpha$ and $\beta$ in the
corresponding interval for $x$.
Another particularly useful result is as follows:

\begin{lem}
\label{lem4}
If a cumulative distribution function $F(x)$ is log-concave
on interval $C = (a, b)$, then $F(x) - F(x_0)$ is also log-concave 
for $\max(a, x_0) < x < b$.
Similarly, $F(x_0) - F(x)$ is log-concave
on $ a < x < \min(x_0, b)$.
\end{lem}

\noindent
{\sc Proof}: Note that
\[
\frac{d\log \{F(x)-F(x_0)\}}{dx} 
=
\frac{f(x)}{F(x) - F(x_0)}
=
\frac{f(x)}{F(x)} \left \{ 1 - \frac{F(x_0)}{F(x)}\right \}^{-1}.
\]
Because $F(x)$ is log-concave over $C$, $f(x)/F(x)$ is a decreasing
function. At the same time,
$1 - F(x_0)/F(x)$ is an increasing function. Hence,
\(
{d\log \{F(x)-F(x_0)\}}/{dx} 
\)
is a decreasing function in $\max(a, x_0) < x < b$.
Consequently, 
$F(x) - F(x_0)$ is log-concave.

The proof of the second conclusion is the same.
\hfill{$\blacksquare$}

\begin{thm}
\label{thm3}
Let $X$ be a random variable with cumulative distribution
function $F(x)$. If $F(x)$ is log-concave on interval $C$,
then $\var \{X | X \in A\} \leq \var \{X | X \in B\}$ for any
intervals $A \subset B \subset C$.
\end{thm}

\noindent
{\sc Proof}: 
Since $F(x)$ is log-concave for $x \in C$, by Lemma \ref{lem4},
so is $F(x) - F(a)$ for $x > a$ and $a \in C$.
Applying Lemma \ref{lem3} twice, we find that
\[
\int_{a < x < y < b} \{F(x)-F(a)\} dx dy
\]
is log-concave in $b$ over $b \in C$ for any given $a \in C$.
That is, Condition (\ref{ConditionA}) is satisfied.
Similarly, Condition (\ref{ConditionB}) is also satisfied.
The result then follows from Theorem \ref{thm2}.
\hfill{$\blacksquare$}

\vspace{2ex}

Theorem 3 is almost a special case of Mailhot (1987)
except for a difference in conditions:
Mailhot (1987) requires the density function $f(x)$ rather than 
the cumulative distribution function $F(x)$ 
to be log-concave.
As pointed out in Lemma \ref{thm3},  when the density function $f(x)$
is log-concave so is $F(x)$.
A large number of well-known distributions have a log-concave
density function. In particular, the normal distribution with
any mean and variance is log-concave. Hence,
the normal result in Chen, van Eeden, and Zidek (2010) is a special case
of Mailhot (1987), and also of Theorem 3.

Many commonly used distributions have log-concave density
or log-concave cumulative distribution functions.
We selectively point out that normal, logistic, double exponential,
Weibull ($c x^{c-1}\exp (- x^c)$) and Gamma ($x^{c-1}\exp( - x)$) with $c > 1$
have log-concave densities. Log-normal, Weibull, and Gamma with
$0 < c < 1$ have log-concave cumulative distribution functions; see 
Bagnoli and Bergstrom (2005) for a more complete list.
In short, the sufficient condition of the above theorem is broadly applicable.

\subsection{Cauchy distribution and symmetric distributions}
\label{sect.cauchy}
Let $X$ be a random variable with a standard Cauchy distribution. 
Its cumulative distribution function
$F(x) = \pi/2 + \arctan (x)$ is log-concave in $C= [0, \infty)$.
Hence,
\[
\var \{X | X \in A\} \leq \var \{X | X \in B\}
\]
for any finite intervals $A \subset B \subset [0, \infty)$.
This is a particularly interesting example 
because the variance
of the Cauchy distribution does not exist, and
its density function or its cumulative distribution function is not log-concave.
Numerical investigations indicate that Condition
(\ref{ConditionA}) is not satisfied by the Cauchy distribution.
Hence, $\var \{X | X \in A\}$ is not partially monotonic
in general.

Suppose $X$ is a positive random variable with decreasing
density function over $[0, \infty)$. Then its cumulative
distribution function is easily verified to be log-concave. 
Hence, $\var \{X | X \in A\}$ is partially monotonic on $A \subset [0, \infty)$.
In particular, $\var(X| 0 < X < b)$ is an increasing function of $b$.

Let $X$ be a symmetrically and absolutely continuously
distributed random variable.
Let $f_a(x)$ be the conditional density function of $|X|$
given $|X| < a$.
For any $0 < a < b$
\[
\frac{f_b(x)}{f_a(x)} = \left \{
\begin{array}{ll}
\frac{P(|X|<a)}{P(|X|<b)} & 0 < x \leq a;\\
\infty & a < x < b.
\end{array}
\right .
\]
This implies $|X_a|$ is smaller than $|X_b|$ in the likelihood
ratio order and hence also in stochastic order.
See Theorem 1.C.1 in Shaked and Shanthikumar (2007, p.~43).
Consequently, $E\varphi(|X_a|) \leq E\varphi(|X_b|)$ for any
increasing function $\varphi(\cdot)$ on $[0, \infty)$, and
$\var(X_a) \leq \var(X_b)$ is a special case. 
The t-distribution, with Cauchy as a special case, is an example.

\subsection{Discrete distributions}
\label{sect.discrete}
Can Theorem \ref{thm3} be
generalized to discrete distributions?
We study this problem for integer-valued random variables. 
Since integration in a discrete space becomes summation,
there is hope that the approach for continuous random
variables is still applicable.
While the same approach can be used, we find that
the conditions corresponding to (\ref{ConditionA}) and
(\ref{ConditionB}) are too complex to be insightful.
However, a simple albeit less
general sufficient condition can be obtained.

Let $X$ be an integer-valued random variable
and $p(x)$ its probability mass function. That is, 
\(
p(x) = P(X = x)
\)
for all integers $x$. The cumulative distribution function
of $X$ is then $G(x) = \sum_{k \leq x} p(k)$. 
We now define an absolutely continuous random variable $Y$
so that its density function
\begin{equation}
\label{denY}
f(y) = \sum_k p(k) I( k - 0.5 < y < k+0.5).
\end{equation}
In other words, $Y$ is uniform on each interval $(k - 0.5, k+0.5]$
provided $p(k) > 0$.
For any integer $a \leq b$, it is easily verified that
\[
E\{ X| a \leq X \leq b\} = E\{Y| a-0.5 < Y \leq b+0.5\}
\]
and that
\[
E\{ X^2 | a \leq X \leq b\} = E\{Y^2| a-0.5 < Y \leq b+0.5\} - 1/12.
\]
Consequently,
\begin{equation}
\label{link}
\var \{ X| a \leq X \leq b\} = \var \{Y| a-0.5 < Y \leq b+0.5\} - 1/12.
\end{equation}
We hence have the following lemma.

\begin{lem}
\label{LemDis}
Let $X$ and $Y$ be the two random variables defined earlier.
A sufficient condition for the partial monotonicity of
$\var(X | X \in A)$ is that the
cumulative distribution function of $Y$ satisfies
Conditions (\ref{ConditionA}) and
(\ref{ConditionB}).
\end{lem}

\noindent
{\sc Proof}: If the cumulative distribution function of $Y$ satisfies 
Conditions (\ref{ConditionA}) and (\ref{ConditionB}), then by Theorem \ref{thm3},
the conditional variance $\var(Y| Y \in A)$ is partially monotonic.
Hence, for any finite interval $A$, $\var \{ X| X \in A\}$ is also partially
monotonic by (\ref{link}). This completes the proof.
\hfill{$\blacksquare$}

\vspace{2ex}

It is likely that this condition is also necessary. Yet the condition
for this seemingly neat result is hard to verify so we do not explore further
in this direction. Instead, we strive to find a few 
simple-to-verify sufficient conditions.

If the probability function of $X$ is unimodal, then its corresponding
$Y$ has a monotonic density function on both sides of the mode.
Hence, $\var \{ X| X \in A\}$ is partially monotonic on either side of the mode.
We have two specific examples for the purpose of illustration.

Suppose $X$ has a geometric distribution.
Since its probability mass function is a decreasing function,
its conditional variance 
$\var \{ X| X \in A\}$ is partially monotonic.

Suppose $X$ has a Poisson distribution with mean $\mu$. 
Then its probability mass function is monotonic for $x > \mu$
and for $x < \mu$. Therefore,
$\var \{ X| \mu < X < b \}$ is an increasing function of $b$
and
$\var \{ X| a < X < \mu \}$ is a decreasing function of $a$.

\section{Partial monotonicity of other measures of uncertainty or information}

In this section, we assume that
$X$ is a random variable with a density function $f(x)$
that is log-concave and differentiable.
Let $X_1$ and $X_2$ be two independent
and identically distributed copies of $X$.
Let $A$ be the event $0 < X_1, X_2 < b$ and
denote $F(x) = P( 0 < X < x)$.
The marginal density function of $U = X_1 - X_2$ given
$A$ is
\[
g(u; b)= \int_{u}^{b}  f(x) f(x - u) dx/F^2(b).
\]

As discussed in the Introduction, the size of $U$
represents the repeatability of an experimental result.
Any increasing function of $|U|$ and its expectation
serves as an index of uncertainty. Thus, it is of
interest to study the properties of $U$ under condition $A$.
We will show that
the log-concavity provides many properties of $X_1 - X_2$,
and partial monotonicity for a number
of uncertainty measures, in particular, the Shannon information.

\begin{lem}
\label{lem4}
The density function $g(u; b)$ is decreasing in $u$ for $u \in [0, b]$
for any $b > 0$.
\end{lem}

\noindent
{\sc Proof}. According to Theorem 1.8 in 
Dharmadhikari and Joag-dev (1988, p.~15), if $X_1$ and $X_2$ are
two independent random variables with the same unimodal distribution,
then $X_1 - X_2$ is also unimodal. Because a distribution with log-concave
density is unimodal, this result leads to the conclusion
of this lemma. The result can also be easily and directly verified by
showing that the derivative of $g(u; b)$ with respect to $u$ is non-positive
for $u > 0$. 
\hfill{$\blacksquare$}

\vspace{2ex}

A bivariate function $K(x, y)$ is totally positive of order 2 (TP$_2$) if, for
every choice of points $(x_1, y_1)$, $(x_2, y_2)$ with $x_1 < x_2$ and
$y_1 < y_2$, we have
\[
K(x_1, y_1)K(x_2, y_2) \geq K(x_1, y_2) K(x_2, y_1).
\]
Let $I(x)$ be an indicator function that is 0 or 1
according to whether $x>0$ or $x \leq 0$.
It is seen that $I(b-u)$ is totally positive.
Being a log-concave density, $f(u-x)I(u-x)$ is also TP$_2$ 
according to
Dharmadhikari and Joag-dev (1988, p.~150). In addition,
by composition formula (2.4) in Karlin (1968, p.~16),
\[
g(u, b) =  \int_0^\infty  I(b-u) \{f(u-x)I(u-x)\} dF(x)/F^2(b)
\]
is also TP$_2$. The following lemma is the simple implication
of the total positivity of $g(u; b)$. For the sake of completeness,
we include a quick proof.
For notational simplicity, we write $g(\dot{u}; b)$ for the partial derivative
of $g$ with respect to $u$. 
\vspace{2ex}

\begin{lem}
For $0 < b_1 < b_2$ and $0 < u < b_1$, we have
\[
g(\dot u; b_1)g(u; b_2) \leq g(u; b_1) g(\dot u; b_2).
\]
\end{lem}

\noindent
{\sc Proof}. Because $g(u; b)$ is TP$_2$, we have for any $\delta > 0$,
\[
g(u; b_1) g(u+\delta; b_2) \geq g(u; b_2) g(u+\delta; b_1).
\]
This implies 
\[
g(u; b_1) \{g(u+\delta; b_2)- g(u; b_2)\} 
\geq 
g(u; b_2) \{ g(u+\delta; b_1) - g(u; b_1)\}.
\]
Dividing both sides by $\delta$ and letting $\delta \to 0$, we get the
result.
\hfill{$\blacksquare$}

\vspace{2ex}
With this result, the following theorem is trivial when coupled with
the concept of the likelihood ratio order introduced earlier.

\begin{thm}
\label{thm4}
For any function $\varphi(u)$ that increases in $|u|$,
\(
E\{ \varphi(U)| X_1, X_2 \in A\}
\)
is partially monotonic.
\end{thm}

\noindent
{\sc Proof}. 
Lemma 5 implies that $\log f(u; b_2) - \log f(u; b_1)$
is an increasing function of $u$ over $u > 0$ for all $0 < b_1 < b_2$.
Hence, $|U|$ has a lower likelihood ratio order given $0 < X_1, X_2 < b_1$
than given $0 < X_1, X_2 < b_2$. This implies the result.
\hfill{$\blacksquare$}
\vspace{2ex}

A special example of $\varphi(\cdot)$ is of particular interest. 
When $\varphi(u) = u^2/2$, we have
$E\{ \varphi(U); A\} = \var(X | X \in A)$.
Hence, we have proved again that $\var(X | X \in A)$
is partially monotonic. It is also easily seen that $E\{|X_1 - X_2| ~|A\}$
is partially monotonic. The quality $E\{|X_1 - X_2|\}$ is known as
a measure of concentration and its corresponding U-statistic
is called Gini's mean difference (Serfling, 1980).

\subsection{Shannon information}

Let $\varphi(u) = - \log g(u; b_2)$ for some $b_2 > 0$.
By Lemma \ref{lem4}, this choice of $\varphi(u)$ is an increasing function 
of $|u|$.
Thus, by Theorem \ref{thm4}, when $b_1 < b_2$ we have
\[
E \{ - \log g(U; b_2)| 0 < X_1, X_2 < b_1 \} 
\leq E \{ - \log g(U; b_2); 0 < X_1, X_2 < b_2 \}.
\]
In other words, 
\begin{equation}
\label{ineq1}
- \int g(u; b_1) \log g(u; b_2)du  \leq - \int g(u; b_2) \log g(u; b_2) du .
\end{equation}
According to Jensen's inequality (Serfling, 1980, p.~351), for any convex function
$\phi(\cdot)$ and random variable $Y$, we have
$E \phi(Y) \geq \phi(E(Y))$ when the expectations exist. 
Applying this inequality to the convex function $- \log (\cdot)$ and random
variable
\[
Y = g(U; b_2)/g(U; b_1)
\]
such that $U$ has density function $g(u; b_1)$,
we get
\[
- \int  g(u; b_1) \log \left [ \frac{g(u; b_2)}{g(u; b_1)} \right ] du
\geq - \log \left [ \int \frac{g(u; b_2)}{g(u; b_1)} g(u; b_1) du \right ]
\geq 0
\]
because the integration of $g(u; b_2)$ over the support of $g(u; b_1)$
is no more than 1.
Hence, we get
\begin{equation}
\label{ineq2}
- \int g(u; b_1) \log g(u; b_1)du \leq - \int g(u; b_1) \log g(u; b_2)du.
\end{equation}
Combining (\ref{ineq1}) and (\ref{ineq2}), we get
\[
- \int g(u; b_1) \log g(u; b_1)du \leq - \int g(u; b_2) \log g(u; b_2)du.
\]
Note that 
\[
E \{ - \log g(U; b_1)| 0 < X_1, X_2 < b_1\}
\]
is the Shannon information of $U$ given $0 < X_1, X_2 < b_1$.
The inequality hence implies that the conditional Shannon information of $U$
increases in $b$.
More formally, we have the following theorem without proof.

\begin{thm}
\label{thm5}
The Shannon information of $U$ given $A$ is a partially
increasing function of interval $A$.
\end{thm}

This result matches our intuition well. Shannon information measures the
amount of uncertainty in a distribution. For
random variables with log-concave density,
the uncertainty in the form of $X \in A$ increases
when $A$ increases. This helps to reduce the uncertainty measured
by the conditional variance or to increase the information
measured in terms of the Shannon information.

\section{Discussion}
One cannot help but conjecture that similar results hold for multidimensional
random subjects. This would be an interesting investigation. Our discussion on
integer-valued random variables has been limited. In addition, there exist
many other versions of information/entropy (Harris, 1982). 
Our result on random variables
with a log-concave density is clearly applicable to most of them.
We hope that our results
will stimulate interest in these areas.
\vspace{2ex}

\noindent
{\bf Acknowledgment}. This research is partially supported by the Natural
Sciences and Engineering Research Council of Canada. The author would
like to thank the referee for insightful comments.

\vspace{2ex}

\noindent
{\bf \large References}
\begin{enumerate}

\item
Bagnoli, M. and Bergstrom, T. (2005).
Log-concave probability and its applications.
{\it Economic Theory} {\bf 26}, 455-469.

\item
Burdett, K. (1996). Truncated means and variances.
{\em Economics Letters} {\bf 52}, 263-267.

\item
Chen, J., van Eeden, C., and Zidek, J.V.\ (2010).
Uncertainty and the conditional variance.
{\it Probability \& Statistics Letters} {\bf 80}, 1764-1770.

\item
Dharmadhikari, S. and Joag-dev, K. (1988).
{\em Unimodality, Convexity and Applications}.
Boston: Academic Press, Inc.

\item
Harris, B.\ (1982). Entropy.
In {\em Encyclopedia of Statistical Science}
Volume {\bf 2}.
Eds. S. Kotz and N.L. Johnson. New York: Wiley, 512-516.

\item
Karlin, S. (1968). {\em Total Positivity.} Stanford, CA: Stanford University Press.

\item
Mailhot, L.\ (1987).
Ordre de dispersion et lois tronqu\'{e}es, 
C. R. Acad. Sci. Paris, S\'{e}r. I Math. {\bf 304}, 499-501.

\item
Serfling, R.J.\ (1980). 
{\em Approximation Theorems in Mathematical Statistics}.
New York: John Wiley \& Sons Inc.

\item
Shaked, M.\ and Shanthikumar, J.G.\ (2007).
{\em Stochastic Orders}.
New York: Springer Series in Statistics.

\item
Shannon, C.E.\ (1948).
A mathematical theory of communication.
{\em Bell System Technical Journal} {\bf 27}, 379-423, 623-656.

\item
Zidek, J.V.\ and van Eeden, C. (2003).
Uncertainty, entropy, variance and the effect of partial information.
In {\em Mathematical Statistics and Applications: Festschrift for Constance
van Eeden}. Eds. M. Moore, S. Froda, and C. Leger. 
Lecture notes-monograph series, Institute of Mathematical Statistics, 
Volume {\bf 42}, 155-167.

\end{enumerate}

\end{document}